\documentclass{commat}

\usepackage{caption}
\usepackage{subcaption}

\title{%
    On Young diagrams of maximum dimension
    }

\author{%
    Vasilii Duzhin, Egor Smirnov-Maltsev
    }

\authorinfo[%
    V. Duzhin]{
    Saint Petersburg Electrotechnical University ``LETI'', Russia}{%
    vsduzhin@etu.ru
    }

\authorinfo[%
    E. Smirnov-Maltsev]{
    RUDN University, Russia}{%
    1032212272@rudn.ru
    }

\abstract{%
We study the problem of finding Young diagrams of maximum dimension, i.~e. those with the largest number of Young tableaux of their shapes. Consider a class of Young diagrams that differ from a symmetric diagram by no more than one box $(i,j)$ in each row and column. It is proven that when moving boxes $(i,j), i>j$ to symmetric positions $(j,i)$, the original diagram is transformed into another diagram of the same size, but with a greater or equal dimension. A conjecture is formulated that generalizes the above fact to the case of arbitrary Young diagrams. Based on this conjecture, we developed an algorithm applied to obtain new Young diagrams of sizes up to 42 thousand boxes with large and maximum dimensions.
    }



\VOLUME{31}
\YEAR{2023}
\NUMBER{3}
\firstpage{33}
\DOI{https://doi.org/10.46298/cm.12641}

\begin{document}

\section{Introduction} \label{sec:intro}

One of the important problems of asymptotic representation theory is the problem of finding irreducible representations of the symmetric group with maximum dimensions~\cite{bbrock, mckay76}. Since such representations are parameterized by Young diagrams, the problem is equivalent to searching among Young diagrams of a fixed size for those in which the largest number of Young tableaux can be placed. In other words, the problem comes down to finding Young diagrams of maximum dimension or \textit{maximum Young diagrams}.

In~\cite{verpavl09}, the first 130 Young diagrams with maximum dimensions were found by an exhaustive search. It was also assumed there that maximum diagrams differ from symmetric ones by no more than one box. All such diagrams were found up to level 310 of the Young graph. Finding larger maximum diagrams is difficult due to the exponential growth of their number.

Different strategies for searching Young diagrams with large and maximum dimensions were proposed in~\cite{ius15, pdmi15, cte19}.
Particularly, a so-called greedy sequence of Young diagrams was constructed from the root of the Young graph. In this sequence, each subsequent Young diagram is obtained from the previous one by adding a box with the maximum Plancherel probability. It was found that although such a sequence contains diagrams with large dimensions, not all of them are maximum. For example, the 15th diagram of the greedy sequence already does not have the maximum dimension. Some algorithms have been proposed whose idea is to iteratively transform a Young diagram into another diagram of the same size, but presumably of greater dimension. These algorithms were applied to Young diagrams from the greedy sequence, resulting in a sequence whose each diagram has the largest dimension among all known diagrams of the same size.

In~\cite{ius15}, the so-called shaking algorithm was used to find diagrams with maximum dimensions. The idea of this algorithm is to add $k$ boxes with the maximum Plancherel probabilities to the diagram and then remove $k$ boxes with the minimum Plancherel probabilities. This algorithm was applied for different values of the parameter $k$.

A similar problem was studied for strict Young diagrams~\cite{pdmi15, knots16} parameterizing projective representations of the symmetric group. Using an exhaustive search algorithm, the first 250 strict diagrams with maximum dimensions were found. Along with the shaking algorithm described above, the so-called branches algorithm was also used. It considers not one, but $m$ greedy sequences, which are constructed from diagrams obtained from the original one by the shaking algorithm.

In the work~\cite{cte19}, in addition to the Plancherel greedy sequence of standard Young diagrams, the so-called greedy blow-up sequence was constructed. In this sequence, each new box is added to the diagram in such a way as to bring its profile as close as possible to the Vershik-Kerov-Logan-Shepp limit shape. Also in ~\cite{cte19} a new variant of the shaking algorithm was proposed. Same as in the original algorithm, $k$ boxes are sequentially added to the diagram, and then $k$ other boxes are removed. In this case, each box added or removed is selected from the list of $m$ optimal boxes, i.~e. having $m$ largest or $m$ lowest Plancherel probabilities, respectively.

Note that the sequence constructed using strategies from~\cite{ius15, pdmi15, cte19} contains all known Young diagrams with maximum dimensions. It also turned out that the standard Young diagrams from this sequence have an interesting geometric property. All the boxes that distinguish them from the symmetric diagram of the nearest smaller size are concentrated on one side of the line $y=x$, and in each row and in each column there is no more than one such box. In this paper, we prove a theorem that describes a fairly wide class of diagrams whose dimensions cannot be maximum. A conjecture was formulated that generalizes the theorem to the case of arbitrary Young diagrams. We propose the strategy for the algorithm for finding new Young diagrams with large and maximum dimensions.

\section{Basic definitions} \label{sec:basics}

A Young diagram of size $n$ is a finite set of $n$ boxes arranged in rows of non-increasing length~\cite{fulton}.
The hook of a box $(i,j)$ in a Young diagram is a set of boxes consisting of the box $(i,j)$ itself, the boxes located above it in the same column, and the boxes located to the right of it in the same row. The reverse hook of a box $(i,j)$ includes the box itself, as well as the boxes located below it in the same column and to the left of it in the same row.

A standard Young tableau~\cite{fulton} is a diagram of size $n$ whose boxes contain the first $n$ natural numbers. In this case, these numbers are arranged in ascending order from bottom to top and from left to right. The dimension of a Young diagram is the number of Young tableaux for a given diagram. It can be calculated using the hook length formula~\cite{hook}:

\begin{equation}
\dim(\lambda_n) = \frac{n!}{\prod\limits_{(i,j) \in \lambda_n} h(i,j)},
\label{eq1}
\end{equation}
where $\lambda_n$ is a diagram of size $n$, $h(i,j)$ is the hook length of a box $(i,j)$, i.~e. the number of boxes in a hook $(i,j)$.

Diagrams $\lambda_n$ and $\lambda^{*}_n$ are called \textit{symmetrically conjugate} if one of them can be obtained from the other by reflecting boxes with respect to the line $y=x$. It is obvious that symmetrically conjugate diagrams have the same dimension since each tableau of the shape $\lambda_n$ has a one-to-one correspondence to some tableau of the shape $\lambda^{*}_n$, obtained by reflecting the original tableau with respect to the line $y=x$.

The Young graph is a directed graded graph whose vertices are Young diagrams. The level $n$ of the graph consists of all diagrams of size $n$. Each edge leaves the diagram $\lambda_i$ of the level $i$ and enters the diagram $\lambda_{i+1}$ of the level $i+1$, where $\lambda_{i+1}$ is obtained by adding one box to $\lambda_i$. The root of the Young graph is a Young diagram of size 1.

You can define Markov processes on the Young graph by assigning probabilities to edges so that for each vertex the sum of outgoing probabilities is equal to 1. In such processes, the probability of a path is equal to the product of the probabilities of the edges along this path. The Markov process of Plancherel growth (also known as the Plancherel process) is of particular interest for studying the dimensions of Young diagrams. This is a central process, i.~e. the probabilities of different paths between a fixed pair of vertices are equal. In addition, the dimension of the Young diagram $\lambda_n$ is equal to the Plancherel probability of the path from the root of the graph to the diagram $\lambda_n$ multiplied by $n!$~\cite{hook}.

The dimensions of Young diagrams grow exponentially as the size of the diagrams increases. Therefore, it is convenient to use normalization to study dimensions. The normalized dimension of a Young diagram is calculated using the formula~\cite{verker85}:

\begin{equation}
c(\lambda_n) = \frac{-1}{\sqrt{n}} \ln \frac{\dim(\lambda_n)}{\sqrt{n!}} = -\sqrt{(n-1)!} \cdot \sum\limits_{i=1}^{n-1}\ln(p(\lambda_i \nearrow \lambda_{i+1})),
\label{ndim}
\end{equation}
where $p(\lambda_i \nearrow \lambda_{i+1})$ is the transition probability from the $\lambda_i$ diagram to the $\lambda_{i+1}$ diagram in the Plancherel process. Due to the centrality of the Plancherel process, the value of $c(\lambda_n)$ does not depend on the choice of path from the root of the Young graph to the diagram $\lambda_n$. It is obvious from~\eqref{ndim} that the smaller the normalized dimension, the larger the exact one. Therefore, we can reformulate the problem of finding a Young diagram of maximum dimension as the problem of finding a Young diagram of minimum normalized dimension. This fact will be used to construct an algorithm for searching for diagrams of large dimensions in section~\ref{sec:algorithm}.

\section{Dimensions of asymmetric Young diagrams} \label{sec:theorem}

A \textit{symmetric diagram} is a Young diagram in which all boxes are symmetric with respect to the straight line $y=x$.
\textit{The base subdiagram} $\lambda_{sym}$ of the Young diagram $\lambda$ is the maximal symmetric subdiagram of the diagram $\lambda$.

\paragraph{Theorem.} Let $\lambda$ be a Young diagram, $\lambda_{sym}$ be its base subdiagram, $A$ be the set of boxes of the diagram $\lambda$ that are not included in its base subdiagram. Let us assume that each column and each row of $\lambda$ contains at most one box from $A$. $A_u$ are boxes of the set $A$ which are located above the line $y=x$, $A_d$ are boxes of the set $A$ which are located below this line. Let $\widetilde{A_d}$ be the set of boxes symmetric to the boxes $A_d$ with respect to $y=x$. Then the dimension of the diagram $\lambda' = \lambda_{sym} \cup A_u \cup \widetilde{A_d}$ is greater than the dimension of the diagram $\lambda$.

\begin{proof}

First, we prove the theorem for the case when each of the sets $A_u$ and $A_d$ consists of one box. Let the original diagram be $\lambda = \lambda_{sym} \cup (k,l) \cup (t,s)$.

Let the box $(k,l)$ be located above $y=x$, and $(t,s)$ be located below it. Then $k<l$ and $t>s$. In Fig.~\ref{fig:image10} the diagram $\lambda_{sym}$ consists of white boxes, the diagrams from Fig.~\ref{fig:10sub1} and Fig.~\ref{fig:10sub2} are diagrams $\lambda$ and $\lambda'$, respectively. 
\begin{figure}[ht]
\centering
\begin{subfigure}{.47\textwidth}
  \centering
  \includegraphics[width=1\linewidth]{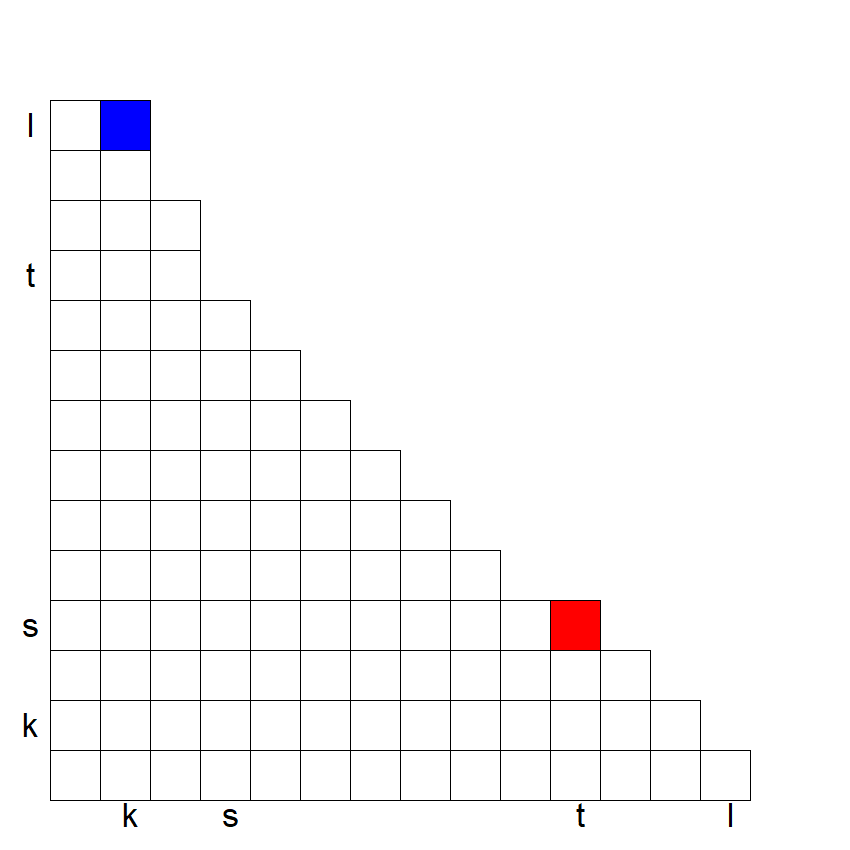}
  \caption{}
  \label{fig:10sub1}
\end{subfigure}%
\begin{subfigure}{.47\textwidth}
  \centering
  \includegraphics[width=1\linewidth]{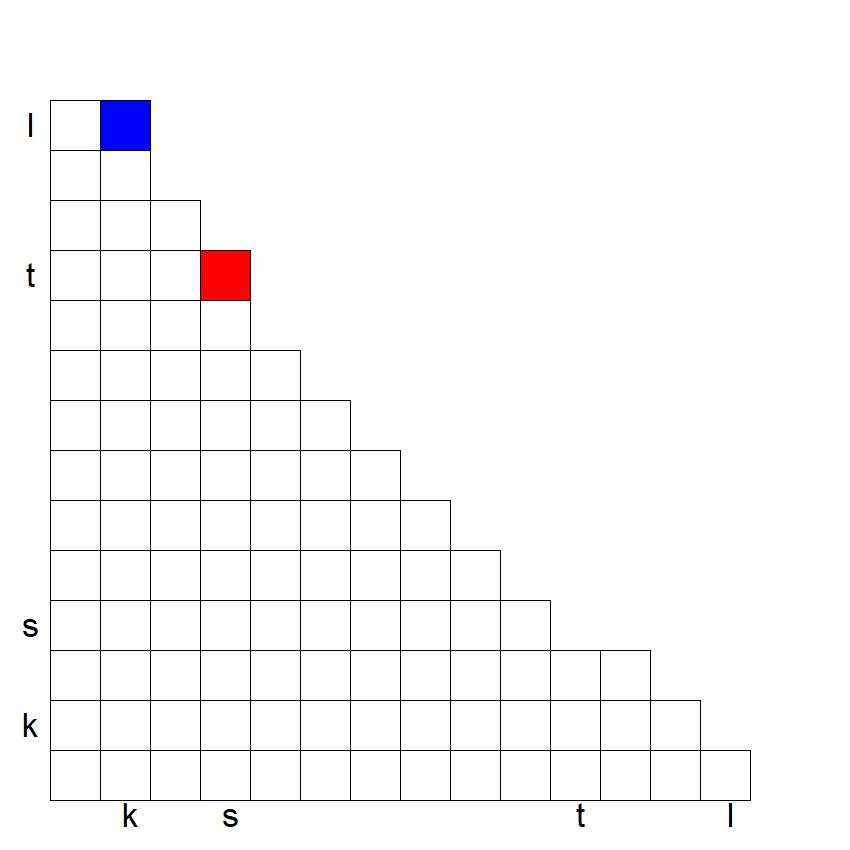}
  \caption{}
  \label{fig:10sub2}
\end{subfigure}
\caption{Transformation of an asymmetric Young diagram to a Young diagram with a larger dimension}
\label{fig:image10}
\end{figure}

Let $h(i,j)$ be the hook length of box $(i,j)$ in the diagram $\lambda$, and let $h'(i,j)$ be the hook length of box $(i,j)$ in the diagram $\lambda'$. Note that $s$ is not equal to $k$, since otherwise the boxes $(k,l)$ and $(t,s)$ would be symmetric to each other with respect to $y=x$, and therefore would be included in the base subdiagram of the $\lambda$ diagram. Similarly, $l$ is not equal to $t$.

Diagrams $\lambda$ and $\lambda'$ consist of the same number of boxes. It follows from the hook length formula~\eqref{eq1} that it is enough to prove that

\begin{equation}
\prod\limits_{(i,j) \in \lambda} h(i,j) > \prod\limits_{(i,j) \in \lambda'} h'(i,j)
\label{eq2}
\end{equation}

Let us prove this inequality. Fig.~\ref{fig:11sub1} and Fig.~\ref{fig:11sub2} show the examples of the original and transformed diagrams, respectively. 
\begin{figure}[ht]
\centering
\begin{subfigure}{.47\textwidth}
  \centering
  \includegraphics[width=1\linewidth]{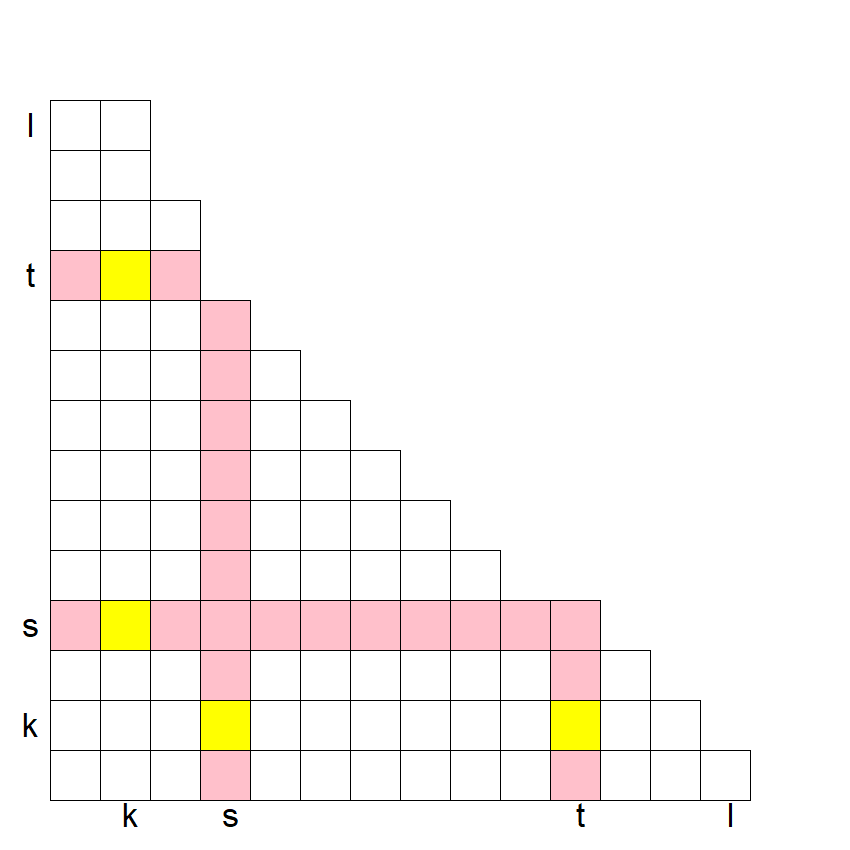}
  \caption{}
  \label{fig:11sub1}
\end{subfigure}%
\begin{subfigure}{.47\textwidth}
  \centering
  \includegraphics[width=1\linewidth]{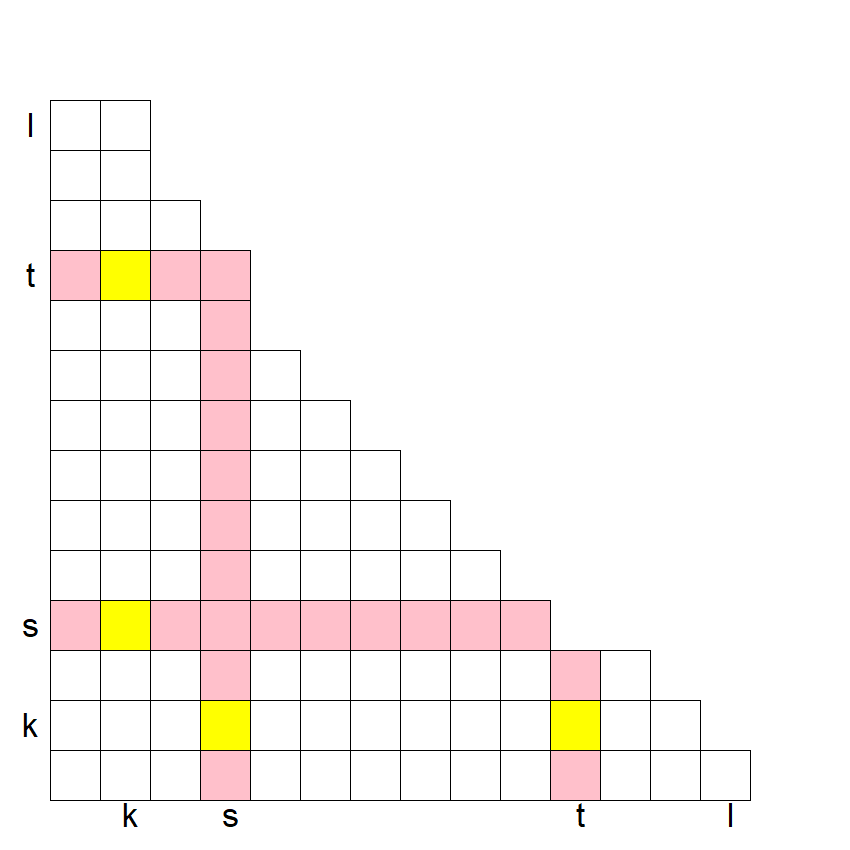}
  \caption{}
  \label{fig:11sub2}
\end{subfigure}
\caption{Different types of boxes of Young diagrams for the case $k<s$}
\label{fig:image11}
\end{figure}

Let us consider 3 types of boxes $(i,j)$ that make up both diagrams.
\begin{enumerate}
\item Boxes $(i,j)$ with coordinates $i \notin \{s, t\}$ and $j \notin \{s, t\}$, i.~e. such that are not included in the reverse hooks of the boxes $(s,t)$ and $(t,s)$. These boxes are highlighted in white in Fig.~\ref{fig:image11}. The hooks of such boxes in the diagrams $\lambda$ and $\lambda'$ do not contain the boxes $(s,t)$ and $(t,s)$, therefore $h'(i,j) = h(i,j)$.

\item Boxes $(i,j)$ belonging to the union of the reverse hooks of boxes $(s,t)$ and $(t,s)$, with coordinates $i \notin \{k, l\}$ and $j \notin \{k, l\}$. In Fig.~\ref{fig:image11} such boxes are highlighted in pink. Let $\lambda_1 = \lambda_{sym} \cup (t,s)$, $\lambda_2 = \lambda_{sym} \cup (s,t)$. The diagram $\lambda_1$ is symmetrically conjugate to the diagram $\lambda_2$, since the subdiagram $\lambda_{sym}$ is symmetric, and the box $(t,s)$ is symmetric with respect to the box $(s,t)$. So $h_{\lambda_1}(i,j) = h_{\lambda_2}(j,i)$. On the other hand, $h(i,j) = h_{\lambda_1}(i,j)$ and $h'(j,i) = h_{\lambda_2}(j,i)$, since the hooks of the boxes $(i,j)$ and $(j,i)$ do not contain the box $(k,l)$. Therefore $h(i,j) = h'(j,i)$.

\item The remaining boxes (marked in Fig.~\ref{fig:image11} in yellow) are boxes, one of the coordinates of which belongs to $\{s,t\}$, and the other coordinate belongs to $\{k, l\}$. Note that the boxes $(s,t)$ and $(t,s)$ are pink, which means that the set of yellow boxes in the diagram $\lambda$ is equal to the set of yellow boxes in the diagram $\lambda'$. Let's call this set $Y$. It is a subset of the set 
$$\{(s,k), (s,l), (t,k), (t,l), (k,s), (l,s), (k,t), (l,t)\}.$$

\end{enumerate}

Recall that if a box $(i,j)$ in a diagram $\lambda$ is white, then $h(i,j) = h'(i,j)$, and if a box $(i,j)$ in $ \lambda$ is pink, then $h(i,j) = h'(j,i)$. 
Hence, we can cancel all factors in inequality~(\ref{eq2}) except those related to yellow boxes. Therefore, it is enough to prove that

\begin{equation}
\prod\limits_{(i,j) \in Y} h(i,j) > \prod\limits_{(i,j) \in Y} h'(i,j)
\label{eq3}
\end{equation}

The box $(k,l)$ in the diagram $\lambda'$ can be either to the left of the box $(s,t)$ or to the right of it.
Let us consider the case when $(k,l)$ is to the left of the box $(s,t)$. This means that $k<s$. An example of such a case is shown in Fig.~\ref{fig:image11}.

Note that $(k,l)$ and $(s,t)$ are corner boxes of $\lambda'$. Therefore, it follows from $k<s$ that $l>t$. So $k<s<t<l$. The boxes $(s,l)$ and $(l,t)$ are not included in the diagram $\lambda$, since even the box $(s,t)$ is not included in it. The boxes $(t,l)$ and $(l,s)$ are not included in $\lambda$, since $(t,s)$ is a corner box. Therefore, $Y = \{(s,k), (t,k), (k,s), (k,t)\}$.

The inequality~\eqref{eq3} in this case can be rewritten as follows
\begin{equation}
h(s,k) \cdot h(t,k) \cdot h(k,s) \cdot h(k,t) > h'(s,k) \cdot h'(t,k) \cdot h'(k,s) \cdot h'(k,t).
\label{eq4}
\end{equation}

Since the lengths of rows and columns with numbers $k, s, t, l$ are known, we can express the hook lengths of yellow boxes in both diagrams as follows:

\begin{equation}
\begin{cases}
h(t,k) = l - t + s - k
\\
h(k,t) = s - k + l - t
\\
h'(t,k) = l - t + s - k - 1
\\
h'(k,t) = l - t + s - k + 1
\\
h'(s,k) = l - s + t - k
\\
h'(k,s) = t - k + l - s
\\
h(s,k) = l - s + t - k - 1
\\
h(k,s) = t - k + l - s + 1.
\end{cases}
\label{case1hooks}
\end{equation}

Expressing all hook lengths in ~\eqref{eq4} in terms of $h'(s,k)$ and $h(t,k)$ using \eqref{case1hooks} and expanding the brackets, we obtain the inequality
$$h'(s,k)^2 > h(t,k)^2.$$
It is true because $s<t$.

Now consider the case when the box $(k,l)$ in the diagram $\lambda'$ is located to the right of the box $(s,t)$, i.~e. $k>s$. It is shown in Fig.~\ref{fig:image13}.

\begin{figure}[ht]
\centering
\begin{subfigure}{.47\textwidth}
  \centering
  \includegraphics[width=1\linewidth]{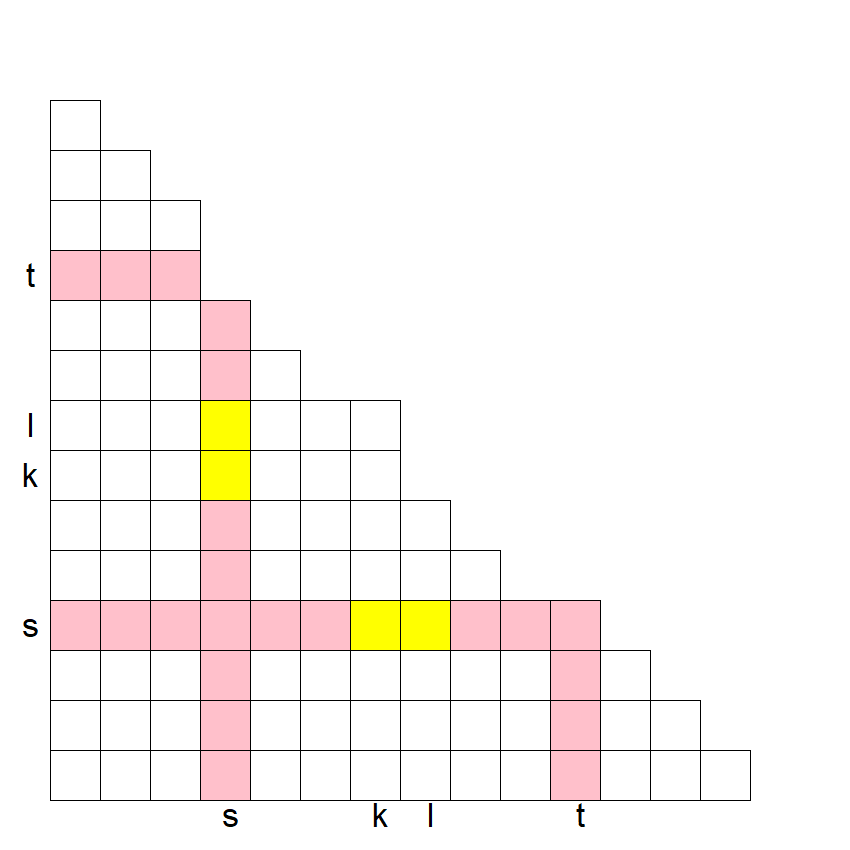}
  \caption{}
  \label{fig:sub1}
\end{subfigure}%
\begin{subfigure}{.47\textwidth}
  \centering
  \includegraphics[width=1\linewidth]{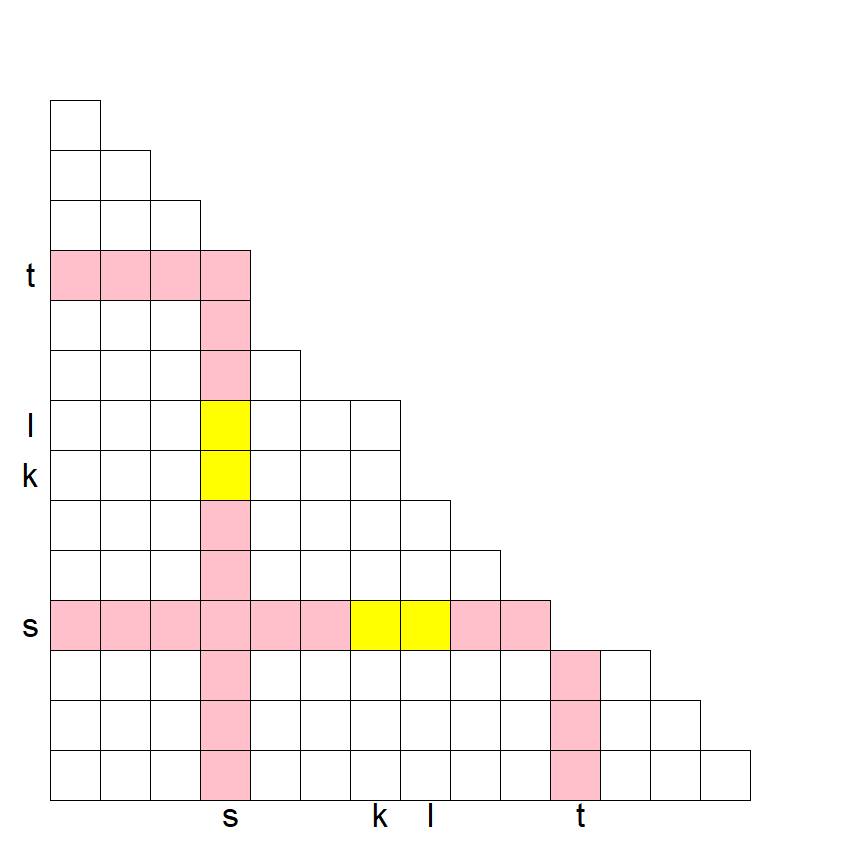}
  \caption{}
  \label{fig:sub2}
\end{subfigure}
\caption{Different types of boxes of Young diagrams for the case $k>s$}
\label{fig:image13}
\end{figure}

The boxes $(k,l)$ and $(s,t)$ are the corner boxes of the diagram $\lambda'$ and $k$ is greater than $s$, therefore $l$ is less than $t$. From this it follows that $s<k<l<t$. Similar to the first case, the boxes $(k,t)$, $(t,k)$, $(t,l)$ and $(l,t)$ are not included in the diagram $\lambda$. Therefore $Y = \{(s,k), (l,s), (k,s), (s,l)\}$.

Inequality~(\ref{eq3}) in this case can be rewritten as follows
\begin{equation}
h(s,k) \cdot  h(l,s) \cdot  h(k,s) \cdot  h(s,l) > h'(s,k) \cdot  h'(l,s) \cdot  h'(k,s) \cdot  h'(s,l).
\label{eq5}
\end{equation}

Let us calculate the hook lengths of yellow boxes in both diagrams using $k,l,s,t$:

\begin{equation}
\begin{cases}
h(s,l) = k - s + t - l
\\
h(l,s) = t - l + k - s
\\
h'(s,l) = k - s + t - l + 1
\\
h'(l,s) = t - l + k - s - 1
\\
h'(s,k) = l - s + t - k
\\
h'(k,s) = t - k + l - s
\\
h(s,k) = l - s + t - k - 1
\\
h(k,s) = t - k + l - s + 1.
\end{cases}
\label{case2hooks}
\end{equation}

Expressing all hook lengths in ~\eqref{eq5} in terms of $h'(s,k)$ and $h(s,l)$ and expanding the brackets, we obtain
$$h'(s,k)^2 > h(s,l)^2.$$
This inequality is true because $k<l$.

Let us proceed to the proof of the theorem in the general case. Let $$A_u = \bigcup\limits_{\alpha}(k_\alpha,l_\alpha), A_d = \bigcup\limits_{\beta}(t_\beta,s_\beta).$$
Then
$$
\begin{cases}
\lambda = \lambda_{sym} \cup \bigcup\limits_{\alpha} (k_\alpha,l_\alpha) \cup \bigcup\limits_{\beta}(t_\beta,s_\beta)
\\
\lambda' = \lambda_{sym} \cup \bigcup\limits_{\alpha} (k_\alpha,l_\alpha) \cup \bigcup\limits_{\beta}(s_\beta,t_\beta).
\end{cases}
$$

In Fig.~\ref{fig:image14} the $\lambda_{sym}$ diagram is made up of white boxes, the diagrams from Fig.~\ref{fig:14sub1} and Fig.~\ref{fig:14sub2} are $\lambda$ and $\lambda'$, respectively.
\begin{figure}[h]
\centering
\begin{subfigure}{.47\textwidth}
  \centering
  \includegraphics[width=1\linewidth]{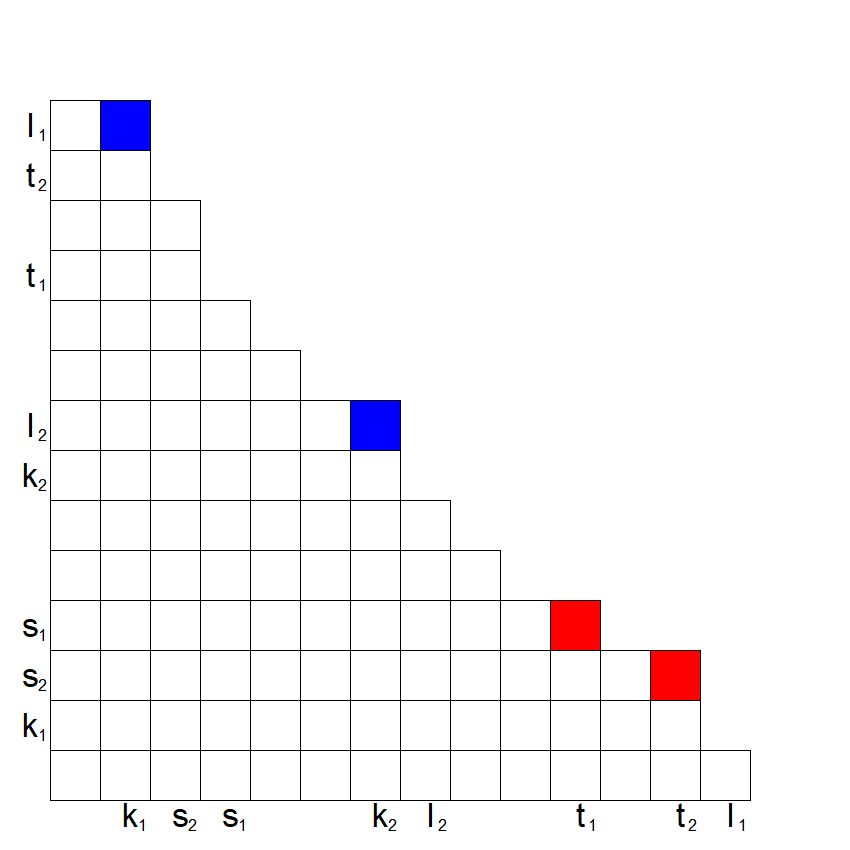}
  \caption{}
  \label{fig:14sub1}
\end{subfigure}%
\begin{subfigure}{.47\textwidth}
  \centering
  \includegraphics[width=1\linewidth]{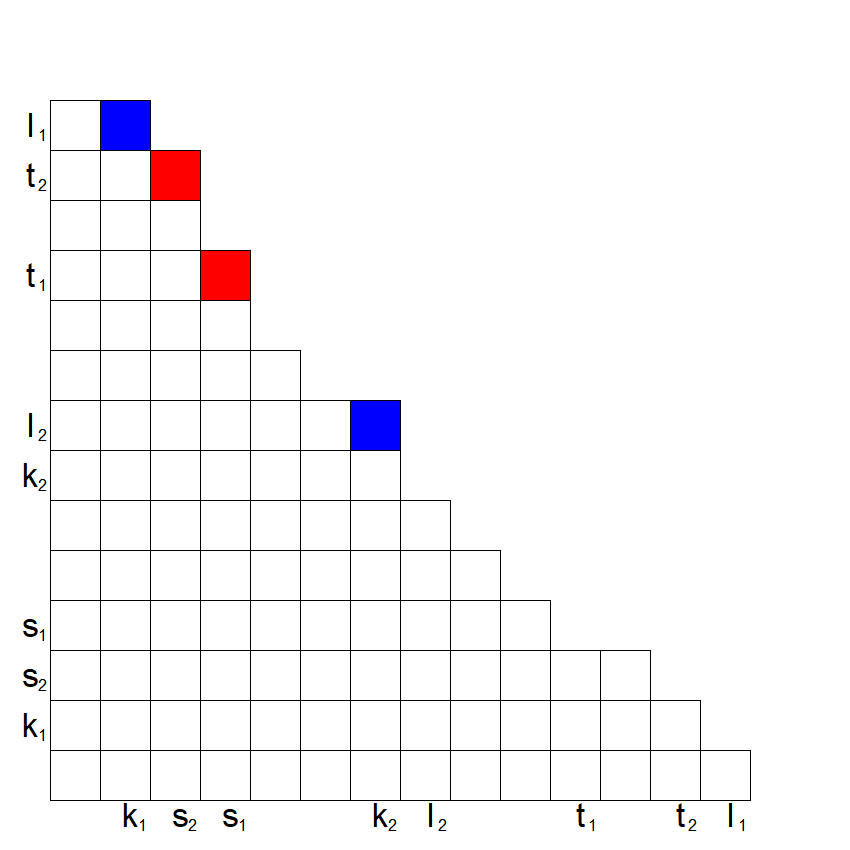}
  \caption{}
  \label{fig:14sub2}
\end{subfigure}
\caption{Transformation of an asymmetric Young diagram to a Young diagram with a larger dimension}
\label{fig:image14}
\end{figure}
Let $h(i,j)$ be the hook length of a box $(i,j)$ in the diagram $\lambda$, and let $h'(i,j)$ be the hook length of this box in $\lambda'$.

By analogy with the case when $A_u$ and $A_d$ consist of one box, we consider three types of boxes of the diagrams $\lambda$ and $\lambda'$. Boxes of different types are highlighted in colours in Fig.~\ref{fig:image15}.

\begin{enumerate}

\item Boxes $(i,j)$ with coordinates $i \notin \bigcup\limits_\beta \{s_\beta, t_\beta\}$ and $j \notin \bigcup\limits_\beta \{s_\beta , t_\beta\}$, i.~e. boxes not lying on the reverse hooks of boxes $(s_\beta,t_\beta)$ and $(t_\beta,s_\beta)$ for any $\beta$. Such boxes $(i,j)$ are highlighted in white in Fig.~\ref{fig:image15}. The hooks of these boxes in the diagrams $\lambda$ and $\lambda'$ do not contain boxes from the sets $\bigcup\limits_\beta(s_\beta,t_\beta)$ and $\bigcup\limits_\beta(t_\beta, s_\beta)$. So $h'(i,j) = h(i,j)$.

\item Boxes $(i,j)$ belonging to the union of the reverse hooks of boxes $(s_\beta,t_\beta)$ and $(t_\beta,s_\beta)$ for any $\beta$, with coordinates $i \notin \bigcup\limits_\alpha \{k_\alpha, l_\alpha\}$ and $j \notin \bigcup\limits_\alpha \{k_\alpha, l_\alpha\}$. In Fig.~\ref{fig:image15} such boxes are highlighted in pink. The hooks of pink boxes do not contain boxes $(k_\alpha,l_\alpha)$ for any $\alpha$. Therefore, the hook length of the pink box $(i,j)$ in the diagram $\lambda$ is equal to the hook length of the box $(j,i)$ in the diagram $\lambda'$.

\item The remaining boxes (marked in yellow) are boxes, one of the coordinates of which belongs to $\bigcup\limits_\alpha \{k_\alpha, l_\alpha\}$, and the other belongs to $\bigcup\limits_\beta \{s_\beta, t_\beta\}$. Note that for any $\beta$ the boxes $(s_\beta,t_\beta)$ and $(t_\beta,s_\beta)$ are pink, which means the set of yellow boxes in the diagram $\lambda$ is equal to the set of yellow boxes in the diagram $\lambda'$. Let us call this set $Y$.
\end{enumerate}
\begin{figure}[ht]
\centering
\begin{subfigure}{.47\textwidth}
  \centering
  \includegraphics[width=1\linewidth]{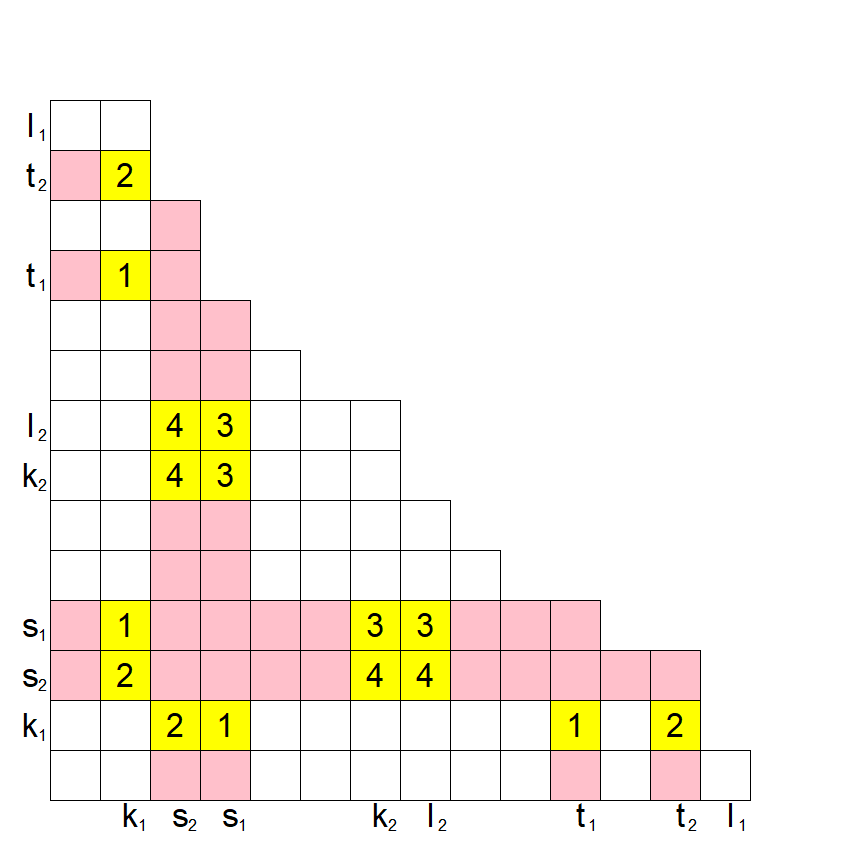}
  \caption{}
  \label{fig:15sub1}
\end{subfigure}%
\begin{subfigure}{.47\textwidth}
  \centering
  \includegraphics[width=1\linewidth]{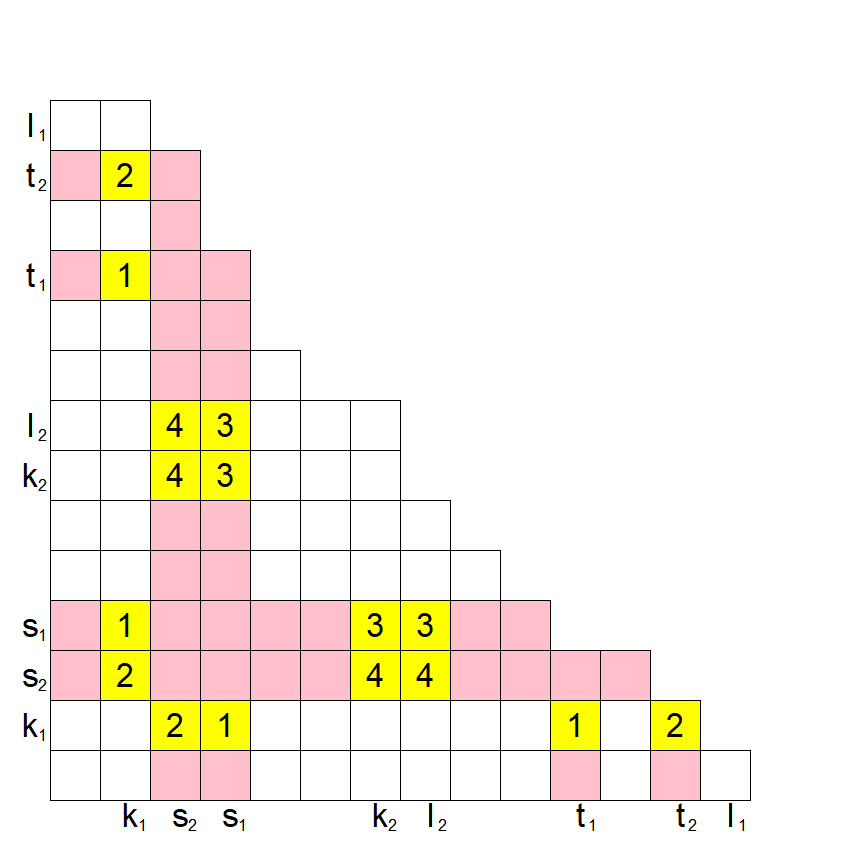}
  \caption{}
  \label{fig:15sub2}
\end{subfigure}
\caption{Different types of boxes of Young diagrams}
\label{fig:image15}
\end{figure}

For each fixed pair of boxes $(\alpha,\beta)$ from $A_u$ and $A_d$, we select a group of four yellow boxes in which one of the coordinates belongs to the set $\{k_\alpha, l_\alpha\}$, and the other belongs to the set $\{s_\beta, t_\beta\}$. Since in the original diagram in each column and in each row no more than one box does not belong to the base subdiagram, the selected groups of boxes for different $\alpha$ and $\beta$ do not intersect. The union of these groups over all pairs $\alpha$ and $\beta$ gives the entire set $Y$. For a fixed pair $(\alpha,\beta)$ we denote such groups as $Y_{\alpha\beta}$. In Fig.~\ref{fig:image15} they are designated by numbers 1, 2, 3, 4.

Since $\lambda$ and $\lambda'$ are the same size and based on the hook length formula~\eqref{eq1} it is enough to prove that
\begin{equation}
\prod\limits_{(\alpha,\beta)} \prod\limits_{(i,j) \in Y_{\alpha\beta}} h(i,j) > \prod\limits_{(\alpha,\beta)} \prod\limits_{(i,j) \in Y_{\alpha\beta}} h'(i,j).
\label{eq6}
\end{equation}
Obviously,~\eqref{eq6} holds if for any pair $(\alpha,\beta)$ the following inequality is true
$$\prod\limits_{(i,j) \in Y_{\alpha\beta}} h(i,j) > \prod\limits_{(i,j) \in Y_{\alpha\beta}} h'(i,j).$$

Consider the diagrams $\lambda_1 = \lambda_{sym} \cup (k_\alpha,l_\alpha) \cup (t_\beta,s_\beta)$ and $\lambda_2 = \lambda_{sym} \cup (k_\alpha ,l_\alpha) \cup (s_\beta,t_\beta)$. For any box $(i,j) \in Y_{\alpha\beta}$ it is true that $h(i,j) = h_{\lambda_1}(i,j)$ and $h'(i,j) = h_{\lambda_2}(i,j)$, since the hooks $h(i,j)$ and $h'(i,j)$ can only contain boxes of the diagram $\lambda_{sym}$, $( k_\alpha,l_\alpha)$, $(t_\beta,s_\beta)$ and $(s_\beta,t_\beta)$, i.~e. only the boxes of the diagrams $\lambda_1$ and $\lambda_2$. So it is enough to prove that
\begin{equation}
\prod\limits_{(i,j) \in Y_{\alpha\beta}} h_{\lambda_2}(i,j) > \prod\limits_{(i,j) \in Y_{\alpha\beta}} h_{\lambda_1}(i,j).
\label{commonw}
\end{equation}
The inequality~\eqref{commonw} corresponds to the case when each of the sets $A_u$ and $A_d$ consists of one box, and was proven earlier.
\end{proof}

In the theorem proved above, the original Young diagram consists of some symmetric subdiagram with at most one box added to each row and each column. We assume that this theorem can be generalized to the case of an arbitrary Young diagram on which no additional conditions are imposed.

\begin{figure}[ht]
\centering
\begin{subfigure}{.5\textwidth}
  \centering
  \includegraphics[width=.6\linewidth]{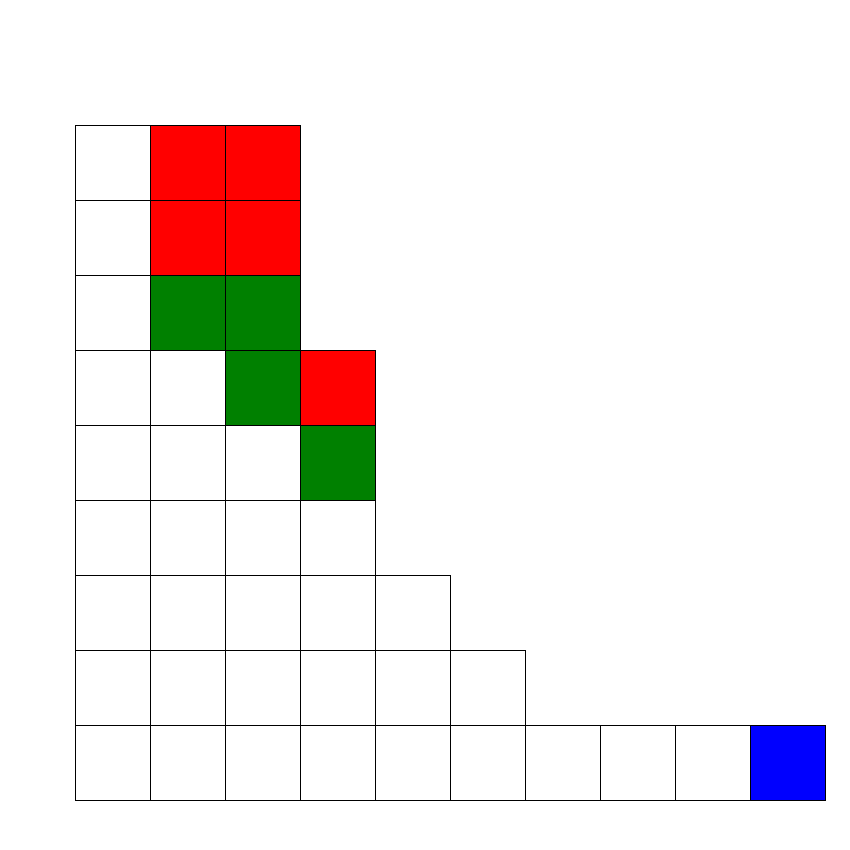}
  \caption{}
  \label{fig:cksub1}
\end{subfigure}%
\begin{subfigure}{.5\textwidth}
  \centering
  \includegraphics[width=.6\linewidth]{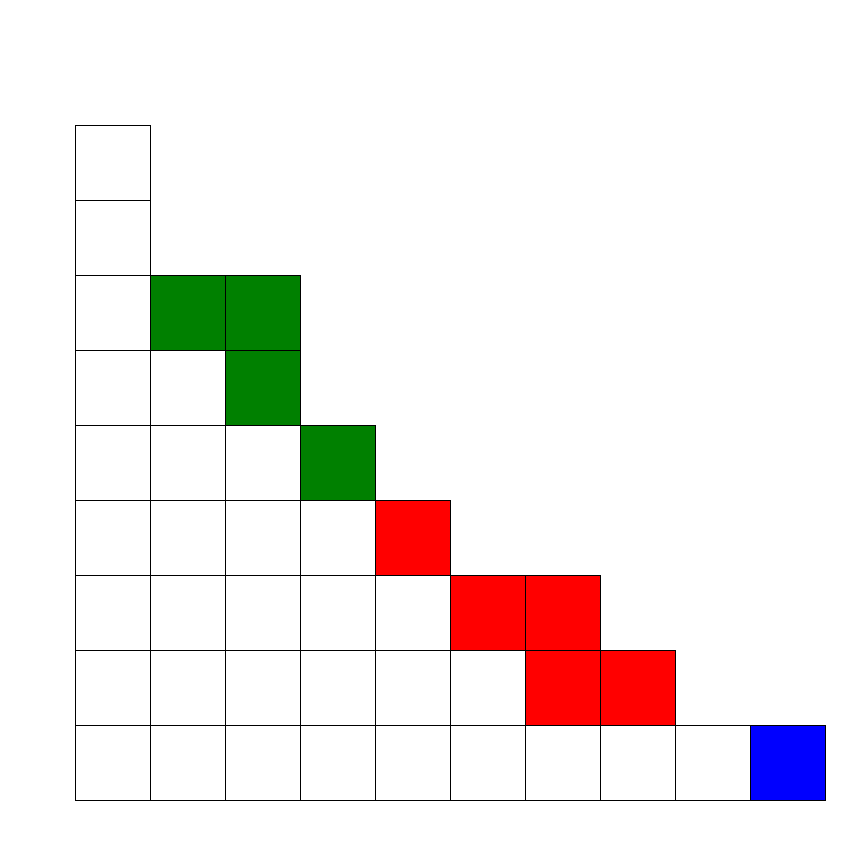}
  \caption{}
  \label{fig:cksub2}
\end{subfigure}
\caption{An example of constructing a Young diagram with a presumably higher dimension}
\label{fig:chet_kl}
\end{figure}

Such a diagram is the union of its base subdiagram $\lambda_{sym}$, boxes $A_u$ lying above the line $y=x$, and boxes $A_d$ lying below this line. Moreover, at least one of the sets $A_u$ and $A_d$ is non-empty. Let $A_u$ be non-empty.
A generalization of the procedure described in the theorem in section~\ref{sec:theorem} is to remove $\lceil \frac{l_i}{2} \rceil$ boxes from the column $i$ containing $l_i$ boxes, followed by the adding the same number of boxes to the row $i$. An example of this transformation is shown in Fig.~\ref{fig:chet_kl}.

Fig.~\ref{fig:cksub1} shows the original diagram $\lambda$, and Fig.~\ref{fig:cksub2} shows the diagram $\lambda'$ obtained after the transformation. The subdiagram $\lambda_{sym}$ is highlighted in white, the boxes $A_u$ that are added and removed are in red, the remaining boxes $A_u$ are in green, and the boxes $A_d$ are in blue. All green boxes, as well as red boxes symmetric to them, are included in the base subdiagram of the $\lambda'$ diagram (see Fig.~\ref{fig:cksub2}). A similar operation can be done if $A_d$ is non-empty, however, boxes will be removed from the row $i$ and added to the column $i$.

The results of the computer experiments described in section~\ref{sec:experiments} give reason to believe that Young diagrams obtained using this transformation always have a larger dimension than the original ones. A consequence of this assumption is the following conjecture about the structure of the Young diagram $\lambda_n$, which has the maximum dimension among all the diagrams of size $n$.

We assume that a Young diagram $\lambda$ of size $n$ of maximum dimension is a disjoint union of its base subdiagram $\lambda_{sym}$ and the set of boxes $A$ lying on one side of the line $y=x$. Moreover, if these boxes lie below the line $y=x$, then each row can contain only one box from $A$. Similarly, if they lie above this line, then each column can contain only one box from $A$.

\section{Algorithm for searching Young diagrams with large and maximum dimensions} \label{sec:algorithm}

The problem of finding Young diagrams of maximum dimension can be formulated as the shortest path problem in the Young graph. This requires to assign a weight $\omega(\lambda_{n-1},\lambda_n) = -\ln(p(\lambda_{n-1} \nearrow \lambda_n))$ to each edge of the graph, where $p(\lambda_{n-1} \nearrow \lambda_n)$ is a Plancherel transition probability of an edge connecting diagrams $\lambda_{n-1}$ and $\lambda_n$. Then the length of any path from the root of the Young graph to the diagram $\lambda_n$ is equal to the normalized dimension~\eqref{ndim} of this diagram divided by $\sqrt{(n-1)!}$. It is known that the smaller the normalized dimension of the Young diagram, the larger its exact dimension. Thus, searching for a diagram of size $n$ of maximum dimension is nothing more than searching for the shortest path from the root of the graph to its level $n$. In order to speed up the algorithm, you can consider not the entire Young graph, but only its subgraph, consisting of diagrams that have certain properties.

Firstly, based on the conjecture formulated in section~\ref{sec:theorem}, it was decided to take into account only diagrams consisting of their base subdiagrams and boxes located below $y=x$. Moreover, each row can contain at most one such box. Consider the subgraph $G'$ of the Young graph consisting exclusively of diagrams satisfying this condition. Obviously, for every diagram $\lambda_n$ in $G'$, there is at least one path from the root to $\lambda_n$. Therefore, we can use search algorithms on $G'$ rather than on the entire Young graph.

Secondly, due to the centrality of the Plancherel process, the costs of all paths from the root of the subgraph $G'$ to the diagram $\lambda_n$ are equal. Therefore, we can consider some spanning tree of the subgraph $G'$. There are different ways to select a spanning tree. We built the so-called \textit{greedy path tree} as follows. Let the diagram $\lambda_n$ be a vertex of $G'$ to which we can add boxes $(i, j)$ and $(i', j')$, and $\dim(\lambda_n \cup (i, j)) > \dim(\lambda_n \cup (i', j'))$. Then we assume that all subtree diagrams rooted at $\lambda_n \cup (i', j')$ will have no box $(i, j)$.

Such a subgraph of the graph $G'$ is indeed a tree since there is a unique path from the root to any vertex $\lambda_n$. This is a path in $G'$ in which at each step that box from $\lambda_n$ is added that maximizes the dimension of the diagram. In such paths, the dimension increases gradually, since the path to the $\lambda_n$ diagram passes exclusively through high-dimensional diagrams.

In this work, the search for the shortest path in this graph is implemented using the $A^*$ algorithm. It uses information about the weights of the graph edges. Namely, for intermediate diagrams $\lambda_m$ the value $f(\lambda_m)=(g(\lambda_m)+h(\lambda_m)) \cdot \sqrt{(n-1)!}$ is calculated, where $g( \lambda_m)$ is the length of the path to the vertex, and $h(\lambda_m)$ is a heuristic estimate of the remaining length. That is, $f(\lambda_m)$ shows the expected normalized dimension of a Young diagram of size $n$ of maximum dimension containing the diagram $\lambda_m$. This value is first calculated from the children of the root. At each next step, the vertex with the minimum $f$ value found is selected and the $f$ function for its descendants is calculated. These steps are repeated until a diagram of size $n$ is reached. This diagram is the result of the algorithm. The absence of a sharp increase in dimensions in the greedy path tree allows us to more accurately determine $h(\lambda_m)$ and thereby reduce the running time of $A^*$.

\section{Computer experiments} \label{sec:experiments}
The $A^*$ algorithm was used to search for Young diagrams of maximum dimension on a tree of greedy paths. In order to estimate the cost of the path from a diagram $\lambda_m$ of size $m$ to level $n$, the function $h(\lambda_m) = P(\lambda_m) \cdot (n-m)$ was calculated, where $P(\lambda_m)$ is the minimum cost of an edge leaving $\lambda_m$ in the greedy path tree. This heuristic does not guarantee that we have found diagrams with maximum dimensions, but it can significantly speed up the process of searching for such diagrams. During the experiments, a sequence of 1000 Young diagrams was constructed. Nine of the resulting diagrams have a larger dimension than diagrams of the same size from the sequence $\Lambda$, which was obtained in the works~\cite{ius15, pdmi15, cte19}.

\begin{figure}[!ht]
	\centering
	\includegraphics[width=.78\textwidth]{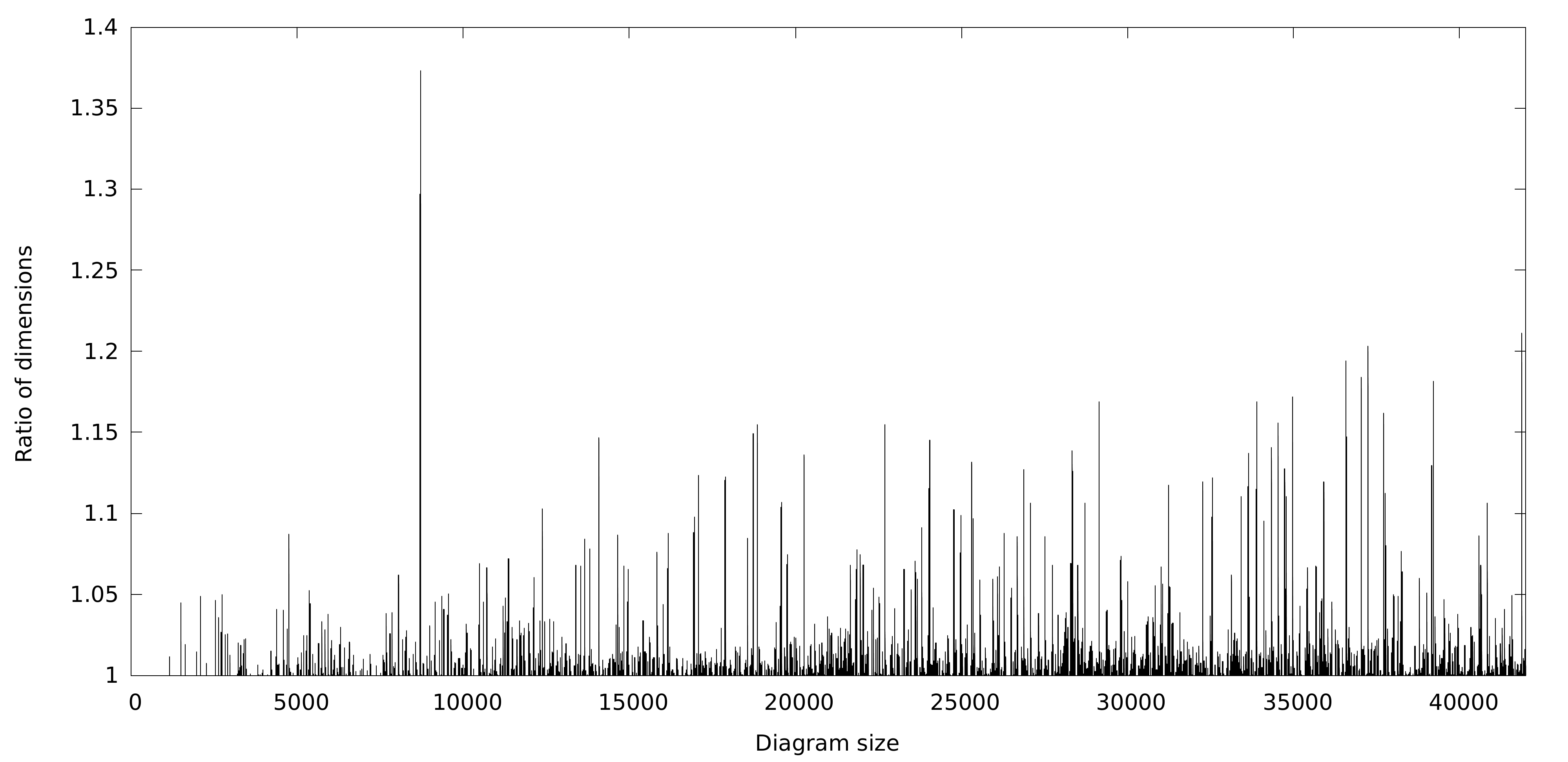}
	\caption{Ratios of dimensions of Young diagrams from the new and original sequences}
	\label{fig:c_ratio}
	\end{figure}

We also searched for the minimum path from some diagram $\lambda_m$ in the greedy path tree to the level $n$. It allows you to find diagrams of maximum dimension that contain the original diagram $\lambda_m$. This approach increases the efficiency of the algorithm but reduces the probability that the found diagram will be the diagram of maximum dimension. Using this strategy, a sequence of Young diagrams was constructed based on the sequence $\Lambda$. Namely, the diagram $\lambda_m$ contained in the original sequence was taken, then the algorithm $A^*$ was run, searching for the minimum path from $\lambda_m$ to level $n=m+3$. Each of the diagrams obtained using $A^*$ was compared with a diagram of the same size from the original sequence. The diagram with the larger dimension was saved. Fig.~\ref{fig:c_ratio} shows the ratio of the exact dimensions of the Young diagrams obtained in this work to the dimensions of the diagrams from $\Lambda$.

In total, in a sequence of length 42 thousand, we managed to find 7278 Young diagrams with dimensions exceeding the dimensions of diagrams of the corresponding size from $\Lambda$.

\section{Conclusion} \label{sec:conclusion}

The article is devoted to the search for Young diagrams with large and maximum dimensions. We consider a class of Young diagrams that differ from a symmetric diagram by no more than one box $(i,j)$ in each row and column.
It is proven that when moving boxes $(i,j), i>j$ to symmetric positions $(j,i)$, the original diagram is transformed into another diagram of the same size, but with a greater or equal dimension. It has been suggested that this fact can be generalized to the case of an arbitrary Young diagram, i.~e. consisting of a symmetric subdiagram with more than one box added to some of its columns and/or rows. A consequence of this assumption is that all boxes of a diagram of maximum dimension that are not included in its maximum symmetric subdiagram are located on one side of the line $y=x$, and each row and column contains at most one such box. This statement has not been strictly proven, but during computer experiments it was found that all known Young diagrams with large and maximum dimensions with sizes up to $10^6$ boxes have this geometric property.

Using this assumption, the search area for diagrams of maximum dimension was reduced. We consider the subgraph $G'$ of the Young graph, consisting exclusively of diagrams that satisfy the above geometric property. It contains a so-called greedy path tree, in which there is only one path from the root of the tree to each vertex. In the subgraph $G'$ this path is greedy, i.~e. at each step the edge with the maximum Plancherel transition probability is selected. The problem of searching for Young diagrams of maximum dimension is reformulated as searching for the shortest path in the resulting tree. This problem solved using the $A^*$ algorithm. As a result, a sequence of 1000 Young diagrams with large dimensions was obtained. Nine of these diagrams have larger dimensions than previously known diagrams of the same size. Also, with the help of $A^*$, an existing sequence of Young diagrams of up to 42 thousand boxes with large dimensions was updated. We found 7278 Young diagrams with dimensions exceeding the dimensions of diagrams of the original sequence.

Future plans include optimizing algorithms for searching for Young diagrams of maximum dimension and adapting them to search for strict Young diagrams of maximum dimension. It is also planned to further study the constructed sequences to identify the properties of Young diagrams of maximum dimension.

\subsection*{Acknowledgments}
This work was supported by grant RSF 22-21-00669. The authors thank prof. Nikolay Vassiliev for numerous discussions and valuable advice, which greatly improved the manuscript.

{\small\bibliography{commat}}
\EditInfo{ December 4, 2023}{ December 5, 2023}{Jacob Mostovoy, Sergei Chmutov}
\end{document}